\documentclass[12pt]{amsart}
\usepackage{amsmath,amssymb,amsbsy,amsfonts,latexsym,amsopn,amstext,cite,
                                               amsxtra,euscript,amscd,bm,mathabx,mathrsfs, abraces}
\usepackage{url}
\usepackage[colorlinks,linkcolor=blue,anchorcolor=blue,citecolor=blue,backref=page]{hyperref}
\usepackage{color}
\usepackage{graphics,epsfig}
\usepackage{graphicx}
\usepackage{float} 
\usepackage[english]{babel}
\usepackage{mathtools}
\usepackage{todonotes}
\usepackage{url}
\usepackage[colorlinks,linkcolor=blue,anchorcolor=blue,citecolor=blue,backref=page]{hyperref}

\usepackage[norefs,nocites]{refcheck}

\hypersetup{breaklinks=true}

\usepackage[norefs,nocites]{refcheck}
\usepackage[english]{babel}
\begin{document}

\newtheorem{thm}{Theorem}
\newtheorem{lem}[thm]{Lemma}
\newtheorem{claim}[thm]{Claim}
\newtheorem{cor}[thm]{Corollary}
\newtheorem{prop}[thm]{Proposition} 
\newtheorem{definition}[thm]{Definition}
\newtheorem{question}[thm]{Open Question}
\newtheorem{conj}[thm]{Conjecture}
\newtheorem{rem}[thm]{Remark}
\newtheorem{prob}{Problem}

\def\ccr#1{\textcolor{red}{#1}}
\def\cco#1{\textcolor{orange}{#1}}

\newtheorem{ass}[thm]{Assumption}

\newtheorem{lemma}[thm]{Lemma}

\newcommand{\GL}{\operatorname{GL}}
\newcommand{\SL}{\operatorname{SL}}
\newcommand{\lcm}{\operatorname{lcm}}
\newcommand{\ord}{\operatorname{ord}}
\newcommand{\Tr}{\operatorname{Tr}}
\newcommand{\Span}{\operatorname{Span}}

\numberwithin{equation}{section}
\numberwithin{thm}{section}
\numberwithin{table}{section}

\def\vol {{\mathrm{vol\,}}}
\def\squareforqed{\hbox{\rlap{$\sqcap$}$\sqcup$}}
\def\qed{\ifmmode\squareforqed\else{\unskip\nobreak\hfil
\penalty50\hskip1em\null\nobreak\hfil\squareforqed
\parfillskip=0pt\finalhyphendemerits=0\endgraf}\fi}

\def \balpha{\bm{\alpha}}
\def \bbeta{\bm{\beta}}
\def \bgamma{\bm{\gamma}}
\def \blambda{\bm{\lambda}}
\def \bchi{\bm{\chi}}
\def \bphi{\bm{\varphi}}
\def \bpsi{\bm{\psi}}
\def \bomega{\bm{\omega}}
\def \btheta{\bm{\vartheta}}
\def \bmu{\bm{\mu}}
\def \bnu{\bm{\nu}}

\newcommand{\bfxi}{{\boldsymbol{\xi}}}
\newcommand{\bfrho}{{\boldsymbol{\rho}}}

\def\vX{\mathbf X}
\def\vY{\mathbf Y}

\def\cA{{\mathcal A}}
\def\cB{{\mathcal B}}
\def\cC{{\mathcal C}}
\def\cD{{\mathcal D}}
\def\cE{{\mathcal E}}
\def\cF{{\mathcal F}}
\def\cG{{\mathcal G}}
\def\cH{{\mathcal H}}
\def\cI{{\mathcal I}}
\def\cJ{{\mathcal J}}
\def\cK{{\mathcal K}}
\def\cL{{\mathcal L}}
\def\cM{{\mathcal M}}
\def\cN{{\mathcal N}}
\def\cO{{\mathcal O}}
\def\cP{{\mathcal P}}
\def\cQ{{\mathcal Q}}
\def\cR{{\mathcal R}}
\def\cS{{\mathcal S}}
\def\cT{{\mathcal T}}
\def\cU{{\mathcal U}}
\def\cV{{\mathcal V}}
\def\cW{{\mathcal W}}
\def\cX{{\mathcal X}}
\def\cY{{\mathcal Y}}
\def\cZ{{\mathcal Z}}
\def\Ker{{\mathrm{Ker}}}
\def\diag{{\mathrm{diag}}}

\def\sA{{\mathscr A}}

\def\NmQR{N(m;Q,R)}
\def\VmQR{\cV(m;Q,R)}

\def\Xm{\cX_m}

\def \A {{\mathbb A}}
\def \B {{\mathbb A}}
\def \C {{\mathbb C}}
\def \F {{\mathbb F}}
\def \G {{\mathbb G}}
\def \L {{\mathbb L}}
\def \K {{\mathbb K}}
\def \N {{\mathbb N}}
\def \Q {{\mathbb Q}}
\def \R {{\mathbb R}}
\def \Z {{\mathbb Z}}

\def \fA{\mathfrak A}
\def \fC{\mathfrak C}

\def \fF{\mathfrak F}
\def \fG{\mathfrak G}
\def \fH{\mathfrak H}
\def \fI{\mathfrak I}
\def \fJ{\mathfrak J}
\def \fL{\mathfrak L}
\def \fM{\mathfrak M}
\def \fP{\mathfrak P}
\def \fR{\mathfrak R}
\def \fS{\mathfrak S}
\def \fZ{\mathfrak Z}

\def \fUg{{\mathfrak U}_{\mathrm{good}}}
\def \fUm{{\mathfrak U}_{\mathrm{med}}}
\def \fV{{\mathfrak V}}
\def \fG{\mathfrak G}
\def \f{\mathfrak G}

\def\e{{\mathbf{\,e}}}
\def\ep{{\mathbf{\,e}}_p}
\def\eq{{\mathbf{\,e}}_q}

 \def\\{\cr}
\def\({\left(}
\def\){\right)}
\def\fl#1{\left\lfloor#1\right\rfloor}
\def\rf#1{\left\lceil#1\right\rceil}

\newcommand{\abs}[1]{\left| #1 \right|}

\def\Im{{\mathrm{Im}}}

\def \ovF {\overline F}

\newcommand{\pfrac}[2]{{\left(\frac{#1}{#2}\right)}}

\def \Prob{{\mathrm {}}}
\newcommand{\disc}{\operatorname{discr}}
\def\e{\mathbf{e}}
\def\ep{{\mathbf{\,e}}_p}
\def\epp{{\mathbf{\,e}}_{p^2}}

\def\Res{\mathrm{Res}}
\def\Orb{\mathrm{Orb}}

\def\vec#1{\mathbf{#1}}
\def \va{\vec{a}}
\def \vb{\vec{b}}
\def \vc{\vec{c}}
\def \vs{\vec{s}}
\def \vu{\vec{u}}
\def \vv{\vec{v}}
\def \vw{\vec{w}}
\def\vlam{\vec{\lambda}}
\def\flp#1{{\left\langle#1\right\rangle}_p}

\def\mand{\qquad\mbox{and}\qquad}

\title{Mean value theorems with smooth numbers}

\author[F.  Majeed]{Fatima Majeed}
\address{Department of Mathematics, Faculty of Computer Science and Ma\-thematics, University of Kufa, 
Najaf Governorate, Iraq}
\email{msc200517@student.uokufa.edu.iq}

\author[I. E. Shparlinski] {Igor E. Shparlinski}
\address{School of Mathematics and Statistics, University of New South Wales, Sydney NSW 2052, Australia}
\email{igor.shparlinski@unsw.edu.au}

  \begin{abstract} We obtain new mean value theorems for exponential sums with very smooth numbers, 
  which provide a power saving against the trivial bound in region where previous bounds do not apply.  \end{abstract}

\subjclass[2020]{11L07, 11N25}

\keywords{smooth numbers, exponential sums, mean values, additive energy}

\maketitle

\tableofcontents

\section{Introduction}
\subsection{Set-up}
\space
A positive integer $n$ is called $y$-smooth if all its prime factors $p$ satisfy $p \le y$.

For $0 \le y \le x$ we consider the set 
\[
\cS(x,y) = \{1 \le n \le x :~ p \mid n, \, \text{prime} \Rightarrow p \le y \}
\]
and as usual denote by 
\[
\Psi(x,y) = \#  \cS(x,y)  
\]
its cardinality. 

Next, for $\vartheta \in [0,1]$ we define the exponential sums
\[
S(\vartheta; x,y) = \sum_{n \in \cS(x,y)} \e(\vartheta n)
\]
where $\e(s) = e^{2\pi i s}$,  
and, for a fixed $\rho > 2$, we consider their mean values 
\[
I_\rho(x,y) = \int_{0}^1 \abs{S(\vartheta; x,y)}^\rho \, d\vartheta
\]
with the goal to improve the trivial bound
\begin{equation}
\label{eq:Triv}
I_\rho(x,y) \le \Psi(x,y)^{\rho-2} I_2(x,y) = \Psi(x,y)^{\rho-1}
\end{equation}
which follows from the Parseval identity 
$ I_2(x,y) = \Psi(x,y) $.

which follows from  the Parseval identity $I_2(x,y)  = \Psi(x, y)$. 

We are especially interested in the {\it additive energy\/}
\[
E(x,y) = I_4(x,y), 
\]
in particular due to the link between the additive energy and 
the Poissonian pair correlation  of the fractional parts, see~\cite[Theorem~1]{ALL}
and a stronger result in~\cite[Theorem~5]{BlWa}. 

Clearly such bounds are only interesting when  $\cS(x,y)$ is a sparse set, 
for example of cardinality of order $x^{1-\kappa}$ for some fixed $\kappa>0$, 
which suggests to consider the values of $y$ of order 
$\(\log x\)^K$ for some fixed $K>0$, see~\cite{Gran,HilTen} about 
the behaviour  of $\Psi(x,y)$. 
In particular, by~\cite[Equation~(1.14)]{Gran}, 
for any fixed $K>1$, we have
\begin{equation}
\label{eq:Growth Psi}
    \Psi(x, \log^K x)=x^{1-\kappa+o(1)}, \qquad  x\to\infty. 
\end{equation}
where throughout the paper we always denote 
\begin{equation}
\label{eq:alpha K}
\kappa = 1/K.
\end{equation}
and for a quantity $U\ge 1$ we write $U^{o(1)}$ for any functions which 
for any fixed $\varepsilon> 0$ is 
bounded by absolute value by $C(\varepsilon) U^\varepsilon$,  where $C(\varepsilon)$ 
depends only on $\varepsilon$.

\subsection{General notation and conventions}  
We recall that  the notations $U = O(V)$, $U \ll V$ and $ V\gg U$  
are equivalent to $|U|\leqslant c V$ for some positive constant $c$, 
which throughout this work, may depend only on the real positive parameter $K$.

We also write $U \asymp  V$ to indicate that  $U \ll V \ll U$. 

We also write $U = V^{o(1)}$ if for any fixed $\varepsilon$ we have

$V^{-\varepsilon} \le |U |\le V^{\varepsilon}$ provided that $V$ is large 
enough. 

For a finite set $\cS$ we use $\# \cS$ to denote its cardinality. 

Throughout this paper, the letter $p$ always denote a prime number.

\subsection{Previous results and new bounds} 

We first recall that a special case of a result of 
Harper~\cite[Theorem~2]{Harp} implies that for some absolute constant $K>0$ and $y \ge \(\log x\)^{K(1 + 1/{\rho-2})}$ we have 
\begin{equation}
\label{eq:HarpBound}
I_\rho(x,y)  \ll \frac{\Psi(x,y)^\rho}{x}, 
\end{equation}
where throughout the paper, the notations $U = O(V)$, $U \ll V$ and $ V\gg U$  are equivalent to $|U|\leqslant c V$ for some positive constant $c$, which may, where obvious,  depends on the real positive parameters $K$ and $\rho$. 

Furthermore,  for an even integer $2s$ the mean value $I_{2s}(x,y)$ 
has an interpretation of the number of solutions to the equation 
\[
n_1+ \ldots+n_s = n_{s+1}+ \ldots+n_{2s}, \qquad n_i\in \cS(x,y), \ 
i =1, \ldots, 2s, 
\]
by a very special case of~\cite[Corollary~16]{BGKS} we have 
\begin{equation}
\label{eq: S-UnitEq}
I_{2s}(x,y) \ll   \Psi(x,y)^{s} +  \Psi(x,y)^{s-1}  \exp\(O\(\frac{y}{\log y}\)\),
\end{equation}
which gives an essentially optimal bound for small $y$. Namely, using a 
classical result of  de Bruijn~\cite[Theorem~1]{dBr} 
(see also~\cite[Theorem~1]{Gran0} or~\cite[Equations~(1.19)]{Gran}),  one can show that
for some absolute constant $c> 0$ and  $y \le c \log x$ the first term 
in~\eqref{eq: S-UnitEq} dominates and thus gives an optimal bound, as obviously 
$I_{2s}(x,y) \gg   \Psi(x,y)^{s}$.

We also recall that  Matthiesen and Wang~\cite{MaWa} have studied more general 
systems of linear equations, however with some additional restrictions and for 
values of $y$ of order at least $\exp\(\sqrt{\log x}\)$, while we are interested the values of 
$y$ that grow as a power of $\log x$. 

Furthermore,  Drappeau and Shao~\cite{DrSha} have considered a Waring-type problem 
with smooth numbers. 

Here we obtain a nontrivial bound in an intermediate range of $y$, 
which interpolates between those of~\cite{BGKS} and~\cite{Harp}. 

It is convenient to define another parameter 
\begin{equation}
\label{eq:beta alpha}
  \beta =  \frac{1-3\kappa}{4} .
\end{equation}

\begin{thm}
\label{thm:MVT-smooth-1} 
Let $y = (\log x)^{K+o(1)}$ for some fixed $K> 3$. 
Then for  $\beta \rho > 1/2$ we have 
\[
 I_\rho(x,y)     \le  x^{o(1)}   \( \Psi(x,y)^{\rho}  x^{-1}   + \Psi(x,y)^{\rho}  x^{-\beta \rho}+  \Psi(x,y)   x^{3(\rho-2)/4}\).
\]  
\end{thm}

\begin{rem}
We have not imposed any conditions on the parameters $K$ and $\rho$
in the formulation of Theorem~\ref{thm:MVT-smooth-1} 
as they are not required in its proof (except for $K> 3$).  However, the region where 
it is really nontrivial is smaller and is given by~\eqref{eq:Nontriv-T1} below. 
\end{rem}

Assuming that $y = (\log x)^{K+o(1)}$,  the first term  in the bound of 
Theorem~\ref{thm:MVT-smooth-1} is always stronger than~\eqref{eq:Triv}, 
and so is the third term, provided $K> 1$. 

For  the second term to be non-trivial, that is, stronger 
 than~\eqref{eq:Triv}, it follows from~\eqref{eq:Growth Psi} that we need 
 \[
\beta \rho >  1- \kappa
 \]
or
\begin{equation}
\label{eq:Large rho}
\rho >  \frac{4(1-\kappa)}{1-3\kappa}. 
\end{equation}

  For the non-triviality of the third term in Theorem~\ref{thm:MVT-smooth-1}  we need
 \[ 
1- \frac{1}{K}  + \frac{3(\rho-2)}{4}   < (\rho -1)\(1- \frac{1}{K}\)
 \]
 or 
 \[
\frac{3}{4}   <  1 -  \kappa.
\]
Thus, recalling also~\eqref{eq:Large rho},  we see that the conditions under which  Theorem~\ref{thm:MVT-smooth-1}
is nontrivial can be summarised as follows:
\begin{equation}
\label{eq:Nontriv-T1}
\rho >  \frac{4(K-1)}{K-3} \mand K > 4.
\end{equation}

In particular,  Theorem~\ref{thm:MVT-smooth-1} does not produce any nontrivial bound on $E(x,y)$ for any  value of $K>1$ and $y = (\log x)^{K+o(1)}$. Thus we now derive a different bound which applies to  $E(x,y)$. 

\begin{thm}
\label{thm:MVT-smooth-2} 
Let $y = (\log x)^{K+o(1)}$ for some fixed $K> 3$. 
Then for  
\[
1/2 < \beta \rho  < 1 \mand  2 < \rho < 2K +4
\]
we have 
\[
 I_\rho(x,y)     \le   \Psi(x,y) ^{1+o(1)}   \(x^{3(\rho-2)/4} + x^{(\rho-2)(1-\eta)} \), 
\] 
where
\[
\eta = \frac{\kappa(\rho-1)}{2 + 3\kappa\rho} .
\]
\end{thm} 

Clearly, Theorem~\ref{thm:MVT-smooth-2} is nontrivial if 
\[
1-\kappa > \max\{3/4, 1-\eta\}
\]
that is, if 
\[
K > 4
\]
and 
\[
 \frac{\rho-1}{2 + 3\kappa\rho}  >1
 \]
which for $\rho>3$ is equivalent to
 \[
 K >   \frac{3\rho}{\rho-3} 
  \]
 (note that  we also need to maintain the condition $K > \rho/2 -2$ in Theorem~\ref{thm:MVT-smooth-2}). 
 
 \begin{cor}
\label{cor:Energy}
Let $y = (\log x)^{K+o(1)}$ for some fixed $K> 12$. 
Then  
\[
E(x,y)     \le  \Psi(x,y) ^{3 - \zeta+o(1)} , 
\]
where
\[
\zeta =  \frac{\kappa(1-12\kappa)}{1 + 5\kappa - 6\kappa^2}. 
\]
\end{cor}
Indeed,  Theorem~\ref{thm:MVT-smooth-2} with the choice $\rho = 4$, becomes
\[
E(x, y) \leq \Psi(x, y)^{1 + o(1)} \( x^{3/2} + x^{2(1 - \eta)} \),
\]
since
\[
2(1 - \eta) = 2\left(1-\frac{3\kappa}{2+12\kappa}\right)=2\left( \frac{2 + 9\kappa}{2 + 12\kappa} \right) = \frac{2+9\kappa}{1+6\kappa} .
\]
 we conclude
\[
E(x, y) \leq \Psi(x, y)^{1 + o(1)} x^{(2+9\kappa)/(1+5\kappa-6 \kappa^2)} . 
\]
Now using~\ref{eq:Growth Psi}  implies Corollary~\ref{cor:Energy}. 

\subsection{Motivation and applications} 

As we have mentioned,  the case of $\rho=4$, that is, bounding $E(x,y)$ is related to the 
so-called metric Poissonian property of smooth numbers, see~\cite{ALL,BlWa}. 
In particular, the bound~\eqref{eq:HarpBound} implies this property for the set of integers $n\ge 2$ 
which are $(\log n)^K$ smooth for some absolute constant  $K>0$.  
Furthermore, using~\cite[Theorem~1.1]{BMS} one can also use a bound on the  
energy $E(x,y)$ to derive a large sieve inequality with smooth moduli.

Unfortunately,  a more explicit bound of Theorem~\ref{thm:MVT-smooth-1}  does not apply to $\rho=4$.

One can also link bounds on $I_\rho(x,y)$ to get 
bounds on the number of solutions to the equation
\[
a_1n_1+ \ldots+a_kn_t =  0 , \qquad n_i\in \cS(x,y), \ 
i =1, \ldots,t, 
\]
uniformly with respect to the non-zero integer coefficients $a_1,  \ldots, a_t$, 
see~\cite[Section~4.1]{MOS}. Then, mimicking  the argument 
of~\cite[Section~5.1]{MOS}, one can  estimate the number of $k\times \ell$ matrices 
of rank $r$ and with their entries from $\cS(x,y)$.

In turn, following the strategy of  Blomer and Li~\cite{BlLi}, one can use such estimates  to  study correlations between the values of linear forms in smooth numbers.

\section{Some known results on  smooth numbers} 

\subsection{Preliminaries}
We first observe that many estimates in the theory of smooth numbers are given in 
term of a certain function $\alpha(x,y)$, known as the {\it saddle point\/} corresponding to the $y$-smooth numbers
up to $x$ as discussed in~\cite{Harp, HilTen}. 
In particular,  $\alpha(x,y)$ satisfies 
\[
\alpha(x,y) = (1+o(1))\frac{\log(1+y/\log x)}{\log y}
\]
provided that $y \leq x$ and $y \to \infty$,  
see~\cite[Theorem~2]{HilTen}.

  We do not need its exact 
definition, but only 
observe that for  $y = (\log x)^{K+o(1)}$, by~\cite[Theorem~2]{HilTen},  we have
\[
\alpha(x,y) = 1 - \kappa + o(1),
\]
see also~\cite[Section~3.10]{Gran} or~\cite[Equation~(2.1)]{Harp}.  

In particular, we often formulate the results below in terms of $\kappa$ rather than 
of $\alpha(x,y)$. However for $y$ of order exacly $\log x$ we need more precise
results in terns of $\alpha(x,y)$

\subsection{Pointwise bounds on exponential sums with smooth numbers} 

We start with recalling a result of Fouvry and Tenenbaum~\cite[Theorem~13]{FoTe}, which we present in a simplified 
form which  suppresses logarithmic factors. 

In what follows,  for $\vartheta \in [0,1]$, $x \ge 1$ 
and integers $a$ and $q \ge 1$, with $\gcd(a,q)=1$
it is convenient to define
\begin{equation}
\label{eq: L gamma}
\cL = 1 + x\abs{\vartheta - \frac{a}{q}} \mand \gamma  = 2 \beta = \frac{1-3\kappa}{2}, 
\end{equation}
where $\kappa$ is given by~\eqref{eq:alpha K} and $\beta$ is given by~\eqref{eq:beta alpha}.  

\begin{lemma}
\label{lem:FT} Let $3 \le y = x^{o(1)}$. For any real $\vartheta$ 
and integers $a$ and $q \ge 1$, with $\gcd(a,q)=1$
we have 
\[
 \abs{S(\vartheta; x,y) } \le x^{1+o(1)}\(x^{-1/4} + q^{-1/2} + (q/x)^{1/2}\) \cL
 \]
as $x \to \infty$. 
\end{lemma} 

We note that the saving here is against $x$ rather than $\Psi(x,y)$. However, in some ranges of $x$ and $y$ such a bound is provided by 
Harper~\cite[Theorem~1]{Harp} (which again we present in a simplified
form.

\begin{lemma}
\label{lem:Harp} Let $y = (\log x)^{K+o(1)}$ for some fixed $K> 3$ and let $\varepsilon > 0$ 
be fixed.  For any real $\vartheta$ 
and integers $a$ and $q \ge 1$, with $\gcd(a,q)=1$ and 
such that 
\begin{equation}
\label{eq:cond theta q}
q \cL \le  2x^{1/2 - \varepsilon}
\end{equation}
we have 
\[
 \abs{S(\vartheta; x,y) } \le \Psi(x,y)   \(q\cL\)^{-\gamma + o(1)}  (\log x)^{5/2+ o(1)} 
  \]
as $x \to \infty$, where $\cL$ and $\gamma$  are given by~\eqref{eq: L gamma}. 
\end{lemma}

A  modified version of Lemma~\ref{lem:Harp} has been given by Baker~\cite[Theorem~2]{Bak}, which 
does not impose any conditions similar to~\eqref{eq:cond theta q} at the cost on an 
additional term in the bound.  Namely, in a simplified form~\cite[Theorem~2]{Bak} yields (note that 
$\kappa$ in~\cite[Theorem~2]{Bak} corresponds to $1-\kappa$ in our notations).

 \begin{lemma}
\label{lem:Bak} Let $y = (\log x)^{K+o(1)}$ for some fixed $K> 3$ and let $\varepsilon > 0$ 
be fixed.  For any real $\vartheta$ 
and integers $a$ and $q \ge 1$, with $\gcd(a,q)=1$ we have 
\[
 \abs{S(\vartheta; x,y) } \le \Psi(x,y)^{1+o(1)}
 \(  \(q \cL \)^{-\gamma} +   \(q \cL  x^{-1+\kappa}\)^{1/2}\)
 \]
as $x \to \infty$,  where $\cL$ and $\gamma$   are given by~\eqref{eq: L gamma}. 
\end{lemma} 

We also notice that there is a series of other bounds~\cite{BrWo,dlB, dlBGr, dlBTen, Drap, DrShp}, however they do not seem to be on help to us.

\section{Proof of Theorem~\ref{thm:MVT-smooth-1}}

\subsection{Initial  splitting}

We fix some parameter $Q > 1$ and define the  sets
\[
\fM_{a,q,Q} =\left \{\vartheta \in [0,1]:~ \abs{\vartheta - \frac{a}{q}} \le \frac{1}{qQ}\right \}.
\]
with $1\le q \le Q$, $1 \le a \le q$ and
$\gcd(a,q)=1$. According to the Dirichlet approximation theorem the sets $\fM_{a,q,Q}$ cover 
the whole interval $[0,1]$. 
Hence
\begin{equation}
\label{eq:I-Maq}
I_\rho(x,y) \le \sum_{q =1}^Q \sum_{\substack{a=1\\\gcd(a,q)=1}}^q 
\int_{\fM_{a,q,Q}}  \abs{S(\vartheta; x,y) }^\rho \, d \vartheta.
\end{equation}

We also note that the Lebesgue measure  $\lambda$  of $\fM_{a,q,Q}$
satisfies 
\[
\lambda\(\fM_{a,q,Q}\) \ll \frac{1}{qQ}.
\]

We now fix some $\varepsilon > 0$ and    split the summation over $q$ in~\eqref{eq:I-Maq}
in the following two sums
\[
I_\rho(x,y)  \le \fI_1(Q) + \fI_2(Q)  , 
\]
where
\begin{align*} 
& \fI_1(Q) =\sum_{1   \le  q  \le x^{1/2-\varepsilon}} \sum_{\substack{a=1\\\gcd(a,q)=1}}^q 
\int_{\fM_{a,q,Q}}  \abs{S(\vartheta; x,y) }^\rho \, d \vartheta,\\
& \fI_2(Q) = \sum_{x^{1/2-\varepsilon} < q \le Q} \sum_{\substack{a=1\\\gcd(a,q)=1}}^q 
\int_{\fM_{a,q,Q}}  \abs{S(\vartheta; x,y) }^\rho \, d \vartheta . 
\end{align*}

{ We assume the condition
\begin{equation}
\label{eq:Large Q}
Q \ge x^{1/2 + \varepsilon}.
\end{equation}

Furthermore, we further observe that for 
 \begin{equation}
\label{eq:small q}
q \le x^{1/2 - \varepsilon} .
\end{equation}
the condition~\eqref{eq:cond theta q}
of  Lemma~\ref{lem:Harp} is satisfied.} Indeed, we have
\[
q\cL = q\left(1 + \frac{x}{qQ} \right) = q + \frac{x}{Q} \le x^{1/2 - \varepsilon} + \frac{x}{x^{1/2 + \varepsilon}} = 2x^{1/2 - \varepsilon}.
\]

\subsection{Bounding  $\fI_1(Q)$} 
\label{sec:bound I1}
 It is convenient to further partition $ \fM_{a,q,Q} $ as
  \begin{equation}
\label{eq:Part M}
\fM_{a,q,Q} = \fM_{a,q,Q}^{\sharp} \bigsqcup\fM_{a,q,Q}^{\flat}
 \end{equation}
with 
\[
\fM_{a,q,Q}^{\sharp}  =\left \{\vartheta \in \fM_{a,q,Q}:~ \abs{\vartheta - \frac{a}{q}} \le \frac{1}{x}\right \}
\]
and 
\[
 \fM_{a,q,Q}^{\flat}= \fM_{a,q,Q} \setminus \fM_{a,q,Q}^{\sharp}. 
\]

We also record the following obvious bounds 
\begin{equation}
\label{eq:Leb-Maq-sharp}
\lambda\(\fM^{\sharp}_{a,q,Q}\)  \ll \frac{1}{x}.
\end{equation}
on  the Lebesgue measure of $\fM^{\sharp}_{a,q,Q}$. 

As we have mentioned, under the conditions~\eqref{eq:Large Q} and~\eqref{eq:small q},
we can apply the bound of Lemma~\ref{lem:Harp}. 
In fact, we consider two cases depending on whether $\vartheta \in \fM_{a,q,Q}^{\sharp} $ or $ \vartheta \in \fM_{a,q,Q}^{\flat} $.

\smallskip
\noindent
 \textbf{Case 1:}  If $ \vartheta \in \fM_{a,q,Q}^{\sharp} $,
then Lemma~\ref{lem:Harp} gives
\[
\abs{S(\vartheta; x,y)} 
\le \Psi(x,y)  (q\cL)^{-\gamma+o(1)}  (\log x)^{5/2 + o(1)} \le  \Psi(x,y) q^{-\gamma} x^ {o(1)},
\]
where we have used the trivial inequality $\cL \ge 1 $.
Then, applying the estimate~\eqref{eq:Growth Psi}, 
we obtain
\begin{equation}
    \label{eq:Bound-S-Harp1}
\abs{S(\vartheta; x,y)} \le \Psi(x,y)^{1 + o(1)} \cdot q^{-\gamma}.
\end{equation}

\smallskip
\noindent
\textbf{Case 2:} If $ \vartheta \in \fM_{a,q,Q}^{\flat} $, we apply Lemma~\ref{lem:Harp},
as before, however we now use that $ \cL \ge x\abs{\vartheta - a/q} $ and obtain
\begin{equation}
    \label{eq:Bound-S-Harp2}
\abs{S(\vartheta; x,y)} \le \Psi(x,y)^{1 + o(1)} \cdot (qx\abs{\vartheta - a/q})^{-\gamma}.
\end{equation}

To estimate $\fI_1(Q)$, we write 
 \begin{equation}
\label{eq:I2-split}
 \fI_1(Q)=\fF^{\sharp}(Q)  +\fF^{\flat} (Q), 
 \end{equation}
 where for $\star \in \{\sharp, \flat\}$ we define 
\[
\fF^{\star}(Q) =\sum_{1\le q \le x^{1/2 - \varepsilon} } \sum_{\substack{a=1\\\gcd(a,q)=1}}^q 
\int_{\fM_{a,q,Q}^\star}  \abs{S(\vartheta; x,y) }^\rho \, d \vartheta. 
 \]

Using~\eqref{eq:Leb-Maq-sharp} ,~\eqref{eq:Bound-S-Harp1} and~\eqref{eq:Bound-S-Harp2} we now derive
\begin{align*}
 \fF^{\sharp}(Q) &   \le   
 \sum_{1\le q \le x^{1/2 - \varepsilon}} \sum_{\substack{a=1\\\gcd(a,q)=1}}^q   \(\Psi(x,y)^{1+o(1)}  q^{-\gamma}\)^\rho
\lambda\(\fM_{a,q,Q}^\sharp\)\\
 &   \le   \Psi(x,y)^{\rho+o(1)}  x^{-1} 
 \sum_{1\le  q \le x^{1/2 - \varepsilon} }  q^{-\gamma \rho+1}\\
 &  \le   \Psi(x,y)^{\rho+o(1)}  x^{-1}\(1+  \(x^{1/2 - \varepsilon} \)^{-\gamma \rho+2}\) \log x.
 \end{align*} 

Absorbing $\log x$ (which actually appears only when $\gamma \rho =  2 $) in $  \Psi(x,y)^{o(1)}$, we obtain
\begin{equation}
\label{eq:Bound-I1-sharp}
 \fF^{\sharp}(Q)    \le    \Psi(x,y)^{\rho+o(1)}  x^{-1} + 
  \Psi(x,y)^{\rho +o(1)}  x^{- \gamma \rho/2+ (\gamma \rho-1)\varepsilon} 
  \end{equation}
  (where we have dropped the factor $\Psi(x,y)^{o(1)} $ in the second term at the cost 
  of increasing the coefficient at the front of $\varepsilon$ in the exponent).

 To bound  $\fF^{\flat}(Q)$, we use ~\eqref{eq:Bound-S-Harp2}, 
 but then we argue slightly differently. 
We write, 
    \begin{align*}
 \fF^{\flat}(Q) &   \le   \Psi(x,y)^{\rho+o(1)}   
 \sum_{1\le  q \le  x^{1/2 - \varepsilon}} \sum_{\substack{a=1\\\gcd(a,q)=1}}^q   (qx)^{- \gamma \rho} \\
 & \qquad \qquad \qquad \qquad \qquad \qquad
\int_{1/x \le  \abs{\vartheta - \frac{a}{q}} \le 1/(qQ)}  \abs{\vartheta - \frac{a}{q}}^{- \gamma \rho} \, d \vartheta \\
&   =\Psi(x,y)^{\rho+o(1)}   
 \sum_{1\le  q \le  x^{1/2 - \varepsilon}} \sum_{\substack{a=1\\\gcd(a,q)=1}}^q   (qx)^{- \gamma \rho} \\ & \qquad \qquad \qquad \qquad \qquad \qquad
\int_{\frac{a}{q} - \frac{1}{qQ}}^{\frac{a}{q} - \frac{1}{x}} 
\left( -\vartheta + \frac{a}{q} \right)^{-\gamma q} \, d\vartheta \\ 
& \qquad \qquad \qquad \qquad \qquad \qquad\qquad \quad
+
\int_{\frac{a}{q} + \frac{1}{x}}^{\frac{a}{q} + \frac{1}{qQ}} 
\left( \vartheta - \frac{a}{q} \right)^{-\gamma q} \, d\vartheta \\
&  =  \Psi(x,y)^{\rho+o(1)}  x^{- \gamma \rho}
 \sum_{1\le  q \le  x^{1/2 - \varepsilon}} \sum_{\substack{a=1\\\gcd(a,q)=1}}^q   q^{- \gamma \rho}  \int_{1/x}^{1/(qQ)} \delta ^{- \gamma \rho} \, d \delta.
 \end{align*}
 Since we have assumed that  $\gamma \rho =  2 \beta \rho > 1$, and we see that 
 \[
  \int_{1/x}^{1/(qQ)} \delta ^{- \gamma \rho} \, d \delta \le x^{\gamma \rho-1}.
  \]
 Thus,  we now derive
\begin{equation}
\begin{split} 
\label{eq:Bound-I1-flat}
 \fF^{\flat}(Q) &   \le   \Psi(x,y)^{\rho+o(1)}  x^{-1}
 \sum_{1\le  q \le  x^{1/2 - \varepsilon}} \sum_{\substack{a=1\\\gcd(a,q)=1}}^q   q^{- \gamma \rho}   \\
&   \le        \Psi(x,y)^{\rho+o(1)}  x^{-1}
 \sum_{1\le  q \le  x^{1/2 - \varepsilon}}    q^{- \gamma \rho+1} \\
  &   \le   \Psi(x,y)^{\rho+o(1)}  x^{-1}\(1+  \(x^{1/2 - \varepsilon} \)^{-\gamma \rho+2}\).
 \end{split}
 \end{equation} 
 Hence we have the same bound on $\fF^{\flat}(Q) $ as on $ \fF^{\sharp}(Q)$ 
 in~\eqref{eq:Bound-I1-sharp}.
Therefore, substituting~\eqref{eq:Bound-I1-sharp}
and~\eqref{eq:Bound-I1-flat} in~\eqref{eq:I2-split} (and dropping $o(1)$ at the cost 
of changing $ (\gamma \rho-12)\varepsilon$ to  $(\gamma \rho-1)\varepsilon$)  we obtain
\begin{equation}
\label{eq:Bound-I1}
\fI_1(Q)  \ll     \Psi(x,y)^{\rho+o(1)}  x^{-1} + 
\Psi(x,y)^{\rho}  x^{- \gamma \rho/2 + (\gamma \rho-1)\varepsilon}. 
 \end{equation} 
 
 Note that there is no dependence on $Q$ in~\eqref{eq:Bound-I1} as long as it
 satisfies~\eqref{eq:Large Q}. 
 
\subsection{Bounding  $\fI_2(Q)$} 
\label{sec:I2}
We notice that for $ x^{1/2-\varepsilon} < q \le Q$, under the assumption~\eqref{eq:Large Q},  we also have $q >x/Q$  
and thus $\fM_{a,q,Q}^{\flat} = \emptyset$.
We also have 
\[
 q^{-1/2} \le   (Q/x)^{1/2} \mand (q/x)^{1/2} \le  (Q/x)^{1/2}.
 \]
In particular, the  bound of Lemma~\ref{lem:FT}  simplifies as 
\[
 \abs{S(\vartheta; x,y) }\ll  x^{1+ o(1)} \(x^{-1/4} + (q/x)^{1/2}\) 
 = x^{ o(1)} \(x^{3/4} + (x Q)^{1/2}\)  . 
\]
  Hence, 
  \begin{align*}
 \fI_2(Q)  
& \ll  x^{ o(1)}   \(x^{3/4} + (xQ)^{1/2}\)^{\rho-2}\
 \sum_{ x^{1/2-\varepsilon}< q \le Q} \sum_{\substack{a=1\\\gcd(a,q)=1}}^q   \\
& \qquad  \qquad  \qquad  \qquad  \qquad  \qquad  \qquad  \quad    \qquad \int_{\fM_{a,q,Q}}  \abs{S(\vartheta; x,y) }^2 \, d \vartheta\\
& \ll  x^{ o(1)}   \(x^{3/4} + (xQ)^{1/2}\)^{\rho-2}
I_2(x,y). 
 \end{align*}
 Using  that  by the Parseval identity $I_2(x,y) = \Psi(x,y)$, we now derive 
  \begin{equation}
\label{eq:Bound-I2}
 \fI_2(Q)  \ll    x^{ o(1)}  \( \Psi(x,y) x^{3\rho/4-3/2} +  
  \Psi(x,y) x^{\rho/2 -1} Q^{\rho/2 -1}\).
 \end{equation}
 
    \subsection{Concluding the proof} 
Combining the bounds~\eqref{eq:Bound-I1}   and~\eqref{eq:Bound-I2}, and dealing in a rather crude manner with terms containing $x^{\varepsilon}$)  we obtain 
 \begin{equation}
\begin{split} 
\label{eq:Bound-I-Q}
 I_\rho(x,y)     \ll    \Psi(x,y)^{\rho+o(1)}  x^{-1} & +\Psi(x,y)^{\rho}   x^{- \gamma \rho/2 + (\gamma \rho-1)\varepsilon} \\
 &  \quad +  \Psi(x,y)^{1+o(1)} x^{3\rho/4-3/2}
 \\
 &   \qquad   
+  \Psi(x,y)^{1+o(1)} x^{\rho/2 -1 } Q^{\rho/2 -1}. 
\end{split}
 \end{equation} 

Clearly  the  bound~\eqref{eq:Bound-I-Q}  is a monotonically 
non-decreasing functions of $Q$. and thus the optimal choice of $Q$ is   
\[
 Q = x^{1/2 + \varepsilon}.
 \]
In this case the  bound~\eqref{eq:Bound-I-Q} becomes 
 \begin{align*}
 I_\rho(x,y)     \ll   \Psi(x,y)^{\rho+o(1)}  x^{-1}  
 +    \Psi(x,y)^{\rho} & x^{-\gamma \rho/2 + ( \gamma \rho -1) \varepsilon} \\
 &  +  \Psi(x,y) x^{3\rho/4 -3/2 +\rho \varepsilon}.
 \end{align*}
Since $\varepsilon>0$  is arbitrary and recalling that $\beta = \gamma/2$, 
we conclude the proof of  Theorem~\ref{thm:MVT-smooth-1}. 
 
\section{Proof of Theorem~\ref{thm:MVT-smooth-2}}

\subsection{Initial  splitting}  We proceed as in the proof of 
Theorem~\ref{thm:MVT-smooth-2}, however we slightly relax~\eqref{eq:Large Q}
as   
 \begin{equation}
\label{eq:Med Q}
Q \ge  x^{1/2} .
\end{equation}
We now fix some $\varepsilon > 0$ and    split the summation over $q$ in~\eqref{eq:I-Maq}
in the following two sums
 
\[
I_\rho(x,y)  \le \fJ_1(Q) + \fJ_2(Q)  ,
\]
where
\begin{align*} 
& \fJ_1(Q) =\sum_{1   \le  q  \le x/Q} \sum_{\substack{a=1\\\gcd(a,q)=1}}^q 
\int_{\fM_{a,q,Q}}  \abs{S(\vartheta; x,y) }^\rho \, d \vartheta,\\
& \fJ_2(Q) = \sum_{x/Q< q \le Q} \sum_{\substack{a=1\\\gcd(a,q)=1}}^q 
\int_{\fM_{a,q,Q}}  \abs{S(\vartheta; x,y) }^\rho \, d \vartheta . 
\end{align*}

For  $\fJ_2(Q)$ we proceed as in Section~\ref{sec:I2} and, since for $Q \ge q > x/Q$, 
we have 
 \[
 \cL \ll 1 \mand  q^{-1/2} + (q/x)^{1/2} \ll (Q/x)^{1/2}, 
 \]
 Lemma~\ref{lem:FT} yields   a full analogue 
of the bound~\eqref{eq:Bound-I2}, that is, 
\begin{equation}
\label{eq:Bound-J2}
\fJ_2(Q)  \ll   x^{o(1)}  \( \Psi(x,y) x^{3\rho/4-3/2} +  
\Psi(x,y) x^{\rho/2 -1} Q^{\rho/2 -1}\).
\end{equation}

\subsection{Bounding  $\fJ_1(Q)$} 
We still partition  $ \fM_{a,q,Q} $ as in~\eqref{eq:Part M}. 

We see from~\eqref{eq:Growth Psi} that Lemma~\ref{lem:Bak}  becomes 
  \begin{equation}
\label{eq:Bound-S-Bak-sharp}
\begin{split}
 \abs{S(\vartheta; x,y) }\le   \Psi(x,y)^{1+o(1)} q^{-\gamma}  + \Psi(x,y)^{1+o(1)}  q^{1/2}  x^{-(1-\kappa)/2}
 \end{split}
  \end{equation}
if $\vartheta \in \fM_{a,q,Q}^{\sharp}$ and 
  \begin{equation}
\label{eq:Bound-S-Bak-flat}
\begin{split}
 \abs{S(\vartheta; x,y) }\le   \Psi(x,y)^{1+o(1)} &q^{-\gamma}  x^{-\gamma}   \abs{\vartheta - \frac{a}{q}}^{-\gamma}
 \\
 &  +  \Psi(x,y)^{1+o(1)}  q^{1/2}  x^{\kappa/2}
  \abs{\vartheta - \frac{a}{q}}^{1/2}
 \end{split}
  \end{equation}
 if $\vartheta \in \fM_{a,q,Q}^{\flat}$. 
 
To estimate $\fJ_1(Q)$, using~\eqref{eq:Leb-Maq-sharp}, \eqref{eq:Bound-S-Bak-sharp} and~\eqref{eq:Bound-S-Bak-flat},
 we write 
 \begin{equation}
\label{eq:J-split}
 \fJ_1(Q) \le \Psi(x,y)^{\rho+o(1)}\(  \fG^{\sharp}(Q) +   \fH^{\sharp}(Q)  + \fG^{\flat} (Q)+ \fH^{\flat}(Q)\), 
 \end{equation}
 where
\begin{align*}
&  \fG^{\sharp}(Q) = 
  \sum_{1\le q \le x/Q } \sum_{\substack{a=1\\\gcd(a,q)=1}}^q  q^{-\gamma \rho} 
  \lambda\(\fM_{a,q,Q}^\sharp\),   \\
 & \fH^{\sharp}(Q) = x^{ - (1-\kappa) \rho /2}
  \sum_{1\le q \le x/Q } \sum_{\substack{a=1\\\gcd(a,q)=1}}^q  q^{ \rho/2} 
\lambda\(\fM_{a,q,Q}^\sharp\) , \\
 &  \fG^{\flat} (Q)=
  x^{-\gamma  \rho}  \sum_{1\le q \le x/Q } \sum_{\substack{a=1\\\gcd(a,q)=1}}^q  q^{-\gamma  \rho} 
\int_{1/x}^{1/(qQ)} \delta ^{-\gamma  \rho} \, d \delta,\\
 &  \fH^{\flat} (Q)=  x^{\kappa \rho /2}
   \sum_{1\le q \le x/Q } \sum_{\substack{a=1\\\gcd(a,q)=1}}^q  q^{ \rho/2} 
\int_{1/x}^{1/(qQ)} \delta ^{ \rho /2} \, d \delta.
\end{align*}
Similarly to the calculations in Section~\ref{sec:bound I1}, and using  out assumption
\[
1\le \gamma \rho = 2 \beta \rho \le 2, 
\] 
we estimate
\begin{align*}
&  \fG^{\sharp}(Q) \ll
  x^{-1} \(1+ (x/Q)^{2-\gamma \rho}\)  
  \ll  x^{1 - \gamma\rho} Q^{-2+\gamma\rho} ,   \\
 & \fH^{\sharp}(Q) \ll
x^{-1  - (1-\kappa)\rho/2} (x/Q)^{2+\rho/2}=   x^{1 + \kappa \rho/2 }  Q^{-2-\rho/2}  ,  
\end{align*}
and
\begin{align*}
 \fG^{\flat} (Q) &\ll  x^{-\gamma  \rho}  \sum_{1\le q \le x/Q } \sum_{\substack{a=1\\\gcd(a,q)=1}}^q  q^{-\gamma  \rho}  x^{\gamma  \rho-1} 
 \ll  x^{-1}  \sum_{1\le q \le x/Q }   q^{-\gamma  \rho+1} \\
& \ll  x^{-1}   \(1+ (x/Q)^{2-\gamma \rho}\)    \ll   x^{1 - \gamma\rho}  Q^{-2+\gamma \rho}. 
\end{align*}
Finally, for  $\fH^{\flat} (Q)$ we obtain
\begin{align*} 
 \fH^{\flat} (Q)&\ll 
 x^{\kappa \rho /2}    \sum_{1\le q \le x/Q } \sum_{\substack{a=1\\\gcd(a,q)=1}}^q  q^{\rho/2} 
(qQ)^{-\rho/2 -1}\\
&=  Q^{-\rho/2 -1}  x^{\kappa \rho /2 + o(1)} \sum_{1\le q \le x/Q }
 \sum_{\substack{a=1\\\gcd(a,q)=1}}^q  q^{-1} \\
& \ll  Q^{-\rho/2 -1}  x^{\kappa \rho /2} (x/Q)^{1+ o(1)} =    x^{1 + \kappa \rho/2 +o(1) }  Q^{-2-\rho/2} .\end{align*}

First we notice that the bounds on $ \fG^{\sharp}(Q)$ and $  \fG^{\flat} (Q) $ coincide
and so do  the bounds on  $ \fG^{\sharp}(Q)$ and $  \fG^{\flat} (Q) $ (up to $x^{o(1)}$. 
Next we compare the bounds on  $ \fG^{\sharp}(Q)$ and  $ \fH^{\sharp}(Q)$. One can easily 
check that 
\[
 x^{1 - \gamma\rho} Q^{-2+\gamma\rho} \ge   x^{1 + \kappa \rho/2 }  Q^{-2-\rho/2} 
 \]
 provided $Q \ge x^{(1-2 \kappa)/(2 - 3 \kappa)}$ which holds under our 
 assumption~\eqref{eq:Med Q}.

Substituting these bounds in~\eqref{eq:J-split} (and noticing that the bounds
on $ \fG^{\sharp}(Q)$ and $  \fG^{\flat} (Q) $ coincide) we now derive
\begin{equation}
\label{eq:Bound-J1}
 \fJ_1(Q)    \le   
    \Psi(x,y)^{\rho+o(1)}  x^{1 - \gamma\rho} Q^{-2+\gamma\rho} .
 \end{equation} 
 
 \subsection{Concluding the proof} 
 We now chose $Q$ to balance (up to factors involving $x^{o(1)}$) the bounds~\eqref{eq:Bound-J2}   and~\eqref{eq:Bound-J1}, that is 
\begin{equation}
\label{eq:Balance}
 Q^{ 1+ \rho(1/2 -\gamma)}\ =    \Psi(x,y)^{\rho-1} x^{2 - (1/2+\gamma)\rho }  .
 \end{equation} 
 It is convenient to  write
 \[
Q = x^\xi
\]
and also observe that 
 \[
 1/2 -\gamma = 3\kappa/2 \mand  1/2+\gamma =1-  3\kappa/2.
\]
Then from~\eqref{eq:Growth Psi} and~\eqref{eq:Balance} we see that an 
asymptotically optimal choice of $\xi$ 
is defined by the equation
\[
\xi (1 + 3\kappa\rho/2) = (1-\kappa) (\rho-1) + 2  - \rho + 3\kappa\rho/2
=  1 +\kappa +\kappa\rho/2. 
\]
Thus, we take 
\[
\xi =  \frac{2 +2\kappa +\kappa\rho}{2 + 3\kappa\rho} = 1-  \frac{2\kappa(\rho-1)}{2 + 3\kappa\rho} = 1- 2\eta.
\]
Obviously we have $\xi< 1$ for $\rho> 1$. On the other hand for $\rho < 2K +4$
we also have $\xi\ge 1/2$, hence for  this choice choice of $\xi$ and $Q=x^\xi$  the condition~\eqref{eq:Med Q} is satisfied and  from~\eqref{eq:Bound-J2} and~\eqref{eq:Bound-J1}
we obtain
\[
  \Psi(x,y)^{\rho}  x^{1 - \gamma\rho} Q^{-2+\gamma\rho}
  = \Psi(x,y) x^{\rho/2 -1} Q^{\rho/2 -1} = \Psi(x,y)  x^{(\rho-2)(1-\eta)+o(1)}, 
 \]
which concludes the proof.

\section*{Acknowledgement}

This work  was  supported, in part,  by the Australian Research Council Grants DP230100530 and DP230100534.


\begin{thebibliography}{99}

\bibitem{ALL} C. Aistleitner, G. Larcher and M. Lewko, `Additive energy and the Hausdorff dimension of the exceptional set in metric pair correlation problems',
\textit{Israel J. Math.},  \textbf{222} (2017), 463--485. 
 
\bibitem{Bak} R. C. Baker, `Smooth numbers in Beatty sequences,
 \textit{Acta Arith.}  \textbf{200} (2021),  429--438.
 

\bibitem{BMS} R. C. Baker, M. Munsch and  I. E. Shparlinski, 
`Additive energy and a large sieve inequality for sparse sequences',
\textit{Mathematika}  \textbf{68} (2022), 362--399. 


 \bibitem{BlLi} V. Blomer and  J. Li, `Correlations of values of random diagonal forms', 
\textit{Intern. Math. Res.  Notices\/} \textbf{2023} (2023),  20296--20336. 

 \bibitem{BlWa} T. F. Bloom and A. Walker, 
 `GCD sums and sum-product estimates', 
\textit{Israel J. Math.}    \textbf{235} (2020), 1--11. 

 \bibitem{BGKS} J.~Bourgain, M.~Z.~Garaev, S. V. Konyagin 
and I. E. Shparlinski, `Multiplicative  congruences with 
variables from short intervals',  
\textit{J. d'Analyse Math.}  \textbf{124} (2014), 117--147.

  \bibitem{BrWo}
  J.~Br\"udern and T.~D Wooley, \emph{Estimates for smooth Weyl sums on major
    arcs}, Intern. Math. Res. Not \textbf{2024} (2024),   14662--14688.

\bibitem{dBr} N. G. de Bruijn, `On the number of positive integers $\leq x$ and free 
of prime factors $>y$, II', 
\textit{ Indag.\ Math.}  \textbf{28} (1966), 239--247.

\bibitem{dlB} R. de la Bret{\`e}che, `Sommes d'exponentielles et entiers sans grand facteur premier',
\textit{Proc. London Math. Soc.}  \textbf{77} (1998), 39--78.


\bibitem{dlBGr} R. de la Bret{\`e}che and A. Granville, `Densit{\'e} des friables',
 \textit{Bull. Soc. Math. France}  \textbf{142} (2014), 303--348.


\bibitem{dlBTen} R. de la Bret{\`e}che and and G. Tenenbaum,  
`Sommes d'exponentielles friables d'arguments rationnels', 
\textit{Funct. Approx. Comment. Math.}  \textbf{37} (2007),  31--38.


 \bibitem{Drap} S. Drappeau, `Sommes friables d'exponentielles et applications',
  \textit{Canad. J. Math.}  \textbf{67} (2015), 597--638.

\bibitem{DrSha} S. Drappeau and X. Shao, 
`Weyl sums, mean value estimates, and Waring's problem with friable numbers',
 \textit{Acta Arith.}  \textbf{176} (2016), 249--299.
 
 \bibitem{DrShp} S. Drappeau and I. E. Shparlinski, 
`Exponential sums  over integers without large prime divisors',
 \textit{Preprint}, 2025, available  from \url{https://arxiv.org/abs/2404.10278}. 

\bibitem{FoTe} {\'E}. Fouvry and G. Tenenbaum, `Entiers sans grand facteur premier en progressions arithmetiques',
\textit{Proc. London Math. Soc.}  \textbf{63} (1991), 449--494.

\bibitem{Gran0}
A.~Granville,  `On positive  integers $\le x$  with prime factors $\le t \log x$', \textit{Number Theory and Applications (Banff, AB, 1988)},   NATO Adv. Sci. Inst. Ser.~C Math. Phys. Sci., v.~265, Kluwer, Dordrecht, 1989, 403--422.

\bibitem{Gran}
A.~Granville, `Smooth numbers: Computational number theory and
beyond', \textit{Proc. MSRI Conf. Algorithmic Number Theory:
Lattices, Number Fields, Curves, and Cryptography, Berkeley
2000}, Cambridge Univ. Press,  267--323. 

\bibitem{Harp}
A. J. Harper, `Minor arcs, mean values, and restriction theory for exponential sums
over smooth numbers',
\textit{Compos. Math.}  \textbf{152} (2016), 1121--1158.

\bibitem{HilTen}
A.~Hildebrand and G.~Tenenbaum,
`On integers free of large prime factors',
\textit{Trans. Amer. Math. Soc.}  \textbf{296}  (1986), 265--290.

 \bibitem{MaWa} L. Matthiesen and M. Wang, 
`Smooth numbers are orthogonal to nilsequences',
\textit{Preprint}, 2023, available  from \url{https://arxiv.org/abs/2211.16892}. 



 \bibitem{MOS} A. Mohammadi, A. Ostafe and I. E. Shparlinski,  `On some matrix counting problems',
\textit{J. Lond. Math. Soc.}  \textbf{110} (2024), Art.~e70044.


\end{thebibliography}
\end{document}